\documentclass[
final
]{dmtcs-episciences}

 \usepackage{amssymb,amsmath,graphicx,graphics,amsthm,amsfonts,mathrsfs,multicol,multirow, amscd,color}
\usepackage[utf8]{inputenc}
\usepackage{subfigure}

\newtheorem{thm}{Theorem}

\newtheorem{lem}[thm]{Lemma}

\newtheorem{con}[thm]{Conjecture}

\usepackage[round]{natbib}

\author{Xiao Zhao\affiliationmark{1}\thanks{Corresponding author. Email: zhaoxiao05@126.com}
  \and Sheng Chen\affiliationmark{2}
}
\title[A note on tight cuts in matching-covered graphs]{A note on tight cuts in matching-covered graphs}
\affiliation{
  College of Mathematics and Information Science, Henan Normal University, Xinxiang, China\\
  Department of Mathematics, Harbin Institute of Technology, Harbin, China}
\received{2020-01-08}
\revised{2020-06-17, 2021-01-19, 2021-05-28}
\accepted{2021-05-28}
\begin{document}
\publicationdetails{23}{2021}{1}{16}{6013}
\maketitle
\begin{abstract}
 Edmonds, Lov\'asz, and Pulleyblank
showed that if a matching covered graph has a nontrivial tight cut,
then it also has a nontrivial ELP-cut. Carvalho et al. gave a stronger
conjecture: if a matching covered graph has a nontrivial tight cut $C$,
then it also has a nontrivial ELP-cut  that does not cross $C$. Chen, et al gave a proof of  the conjecture.
This note is inspired
by the paper of  Carvalho et al. We give a simplified proof of the conjecture, and prove the following result which  is slightly stronger than the conjecture:
if a  nontrivial tight cut $C$ of a matching covered graph $G$ is not an ELP-cut, then
there is a sequence $G_1=G, G_2,\ldots,G_r, r\geq2$  of matching covered graphs,
such that for $i=1, 2,\ldots, r-1$, $G_i$ has an ELP-cut $C_i$, and $G_{i+1}$ is a
$C_i$-contraction of $G_i$, and $C$ is a $2$-separation cut of $G_r$.
\end{abstract}

\section{Introduction}
\label{sec:in}
For graph theoretical terminology and notation,
we  follow \cite{Bondy}. For the terminology that is specific to matching covered graphs,
we follow \cite{ML}. This article studies
finite and undirected loopless graphs.
A {\it perfect matching} of graph $G$ is a set of independent edges which covers all vertices of $G$. An edge of graph $G$ is {\it admissible} if there is a perfect matching of the
graph which contains it. A nontrivial graph is {\it matchable} if it has at least
one perfect matching, and is {\it matching covered} if it is connected and each of
its edges is admissible.  A nontrivial graph $G$ is {\it critical} if $G-v$ is matchable
for any $v\in V(G)$, is {\it bicritical} if $G-u-v$ has a perfect matching
for any two distinct vertices $u$ and $v$ of $G$.

For $X\subseteq V$,  $\overline{X}=V-X$ is the complement of $X$.
The set of all edges
of $G$ with exactly one end in $X$ is denoted by $\partial(X)$, and is referred to as a
cut of $G$. We call $X$
and $\overline{X}$ the shores of $\partial (X)$. For any cut $C:=\partial(X)$ of a graph $G$,
we denote the graph obtained from $G$ by shrinking the shore $X$ to a single
vertex $x$ by $G/(X\rightarrow x)$. Let $G/X$ and $G/\overline{X}$ be obtained from $G$ by contracting $X$
and $\overline{X}$, respectively, and call them $\partial(X)$-contractions of $G$. We call an edge cut $\partial(X)$
trivial if $|X|= 1$ or $|\overline{X}|=1$.  An
edge cut $C:=\partial(X)$ of $G$ is called  a tight cut if $|C\cap M|=1$ for each perfect
matching $M$ of $G$. A matching covered graph which is free of nontrivial tight cuts is a {\it brick} if it is nonbipartite.
% We call a matching covered graph which contains no nontrivial tight cuts is a {\it brace} if it is bipartite, and a {\it brick} otherwise.
If $C:=\partial(X)$ is  a tight cut of a matching covered graph $G$,  then both the $C$-contractions of $G$ are also matching covered. Two cuts $C=\partial(X)$ and $\partial(Y)$ of a matching covered graph cross if each of
the four sets $X\cap Y$, $X\cup \overline{Y}$, $\overline{X}\cap Y$, and $\overline{X}\cap \overline{Y}$ is nonempty. Thus, if $C$ and
$D$ do not cross, then one of the two shores of $C$ is a subset of one of the two
shores of $D$.

A {\it barrier} in a matchable graph $G$ is a subset $B$ of $V$ for which $o(G-B)=|B|$, where $o(G-B)$ is the number of odd components of graph $G-B$. A barrier is trivial if it is a singleton. If $B$ is a barrier and $H$ is any odd connected component of $G-B$, then it is easy to show that $\partial(V(H))$ is tight. Such cuts are
called {\it barrier cuts}. Note that if $B$ is a barrier of matching covered graph $G$, then there is no even component in $G-B$ and $B$ is independent. Two distinct vertices $\{u, v\}$ of matching covered graph $G$ is a $2$-separation if
$G=G_1\cup G_2$ and $V(G_1)\cap V(G_2)=\{u,v\}$ and $|V(G_i)|$ is even for $i=1,2$. It can be verified that both $\partial(V(G_1)-u)$ and $\partial(V(G_2)-v)$ are tight cuts of $G$.  Such cuts are called {\it $2$-separation cuts}. Barrier-cuts and $2$-separation cuts are called {\it ELP-cuts}.

Let $G$ be a matchable graph, and $B$ be a nonempty barrier of $G$.
The bipartite graph $H(B)$ obtained from $G$ by deleting the vertices in the even
components of $G-B$, contracting every odd component to a single vertex,
and deleting the edges with both ends in $B$, is called the {\it core} of $G$ with respect to the barrier $B$. Barrier $B$ of $G$ is a {\it DM-barrier} if each odd component of $G-B$ is critical, and the core $H(B)$ of $G$ with respect to $B$ is matching covered.

\cite{ot} proved the following important result about DM-barrier.

\begin{lem}\textnormal{(\cite{ot})}\label{DM-B}
Let $G$ be a matchable graph, and let $X$ be a nonempty proper subset of $V(G)$ such that both the subgraphs $G[X]$ and $G[\overline{X}]$ are connected, and no edge in the cut $\partial(X)$ is admissible in $G$.
Then $G$ has a DM-barrier $B$ which is a subset of $X$ or of $\overline{X}$. Furthermore, the vertex sets of all the odd components of $G-B$ are also subsets of that same shore.
\end{lem}

\cite{Lov}
proved the ELP Theorem: if a matching covered graph has a nontrivial tight cut,
then it also has a nontrivial ELP-cut. A purely graph
theoretical proof was given by \cite{S}.
\cite{G-e} provided an alternative proof of the ELP Theorem by using Lemma \ref{DM-B}, and gave the following conjecture.

\begin{con}\textnormal{(\cite{ot})}\label{Con}
Let $C$ be a nontrivial tight cut of a matching covered graph
$G$. Then, $G$ has an ELP-cut that does not cross $C$.
\end{con}

\cite{ot} established the validity of  Conjecture \ref{Con} for bicritical graphs
and matching covered graphs with at most two bricks. Also, they observed that Conjecture \ref{Con} is rephrased as follow:
let $C$ be a nontrivial tight cut of a matching covered graph
$G$, then there is a sequence $G_1=G, G_2,\ldots,G_r, r\geq 1$ of matching covered graphs,
such that for $i= 1, 2,\ldots, r-1$, $G_i$ has an ELP-cut $C_i$, $G_{i+1}$ is a
$C_i$-contraction of $G_i$, and  $C$ is an ELP-cut of $G_r$. \cite{Lu} gave a proof of Conjecture \ref{Con} in a preprint.

%\begin{thm}\label{O-component}
%Let $C:=\partial(X)$ be a nontrivial tight cut of matching covered graph $G$. Then, $G$ has a nontrivial $2$-separation cut  does not cross $C$, or has a nontrivial barrier $B$ of $G$  that $B\subsetneq X$ or $B\subsetneq \overline{X}$.
%\end{thm}

Inspired by the proof of Carvalho et al. of the ELP Theorem, we give a simplified proof of Conjecture \ref{Con} and prove a slightly stronger result:

\begin{thm}\label{2-S}
If a  nontrivial tight cut $C$ of a matching covered graph $G$ is not an ELP-cut, then
there is a sequence $G_1=G, G_2,\ldots,G_r, r\geq2$  of matching covered graphs,
such that for $i=1, 2,\ldots, r-1$, $G_i$ has an ELP-cut $C_i$, $G_{i+1}$ is a
$C_i$-contraction of $G_i$, and $C$ is a $2$-separation cut of $G_r$.
\end{thm}

% You may scarsely use \clearpage to advance to a new page if this
% improves the readability of the document structure
%\clearpage
\section{Preliminary}
\label{sec:first}

The following basic results of matching covered graphs will be used  in the proof of our main theorem in the next section. %Carvalho, Lucchesi and Murty\cite{on tight} proved one  particular property of  any nontrivial tight cut  of $G$ with minimal shore.

\begin{lem}\label{2-c}\textnormal{(\cite{2})}
Every matching covered graph is 2-connected except $K_2$.
\end{lem}

\begin{lem}\label{C-t}\textnormal{(\cite{ot})}
Let $G$ be a matching covered graph, and let $C$ be a tight cut of $G$. Then
both $C$-contractions are matching covered. Moreover, if $G'$ is a $C$-contraction
of $G$ then a tight cut of $G'$ is also a tight cut of $G$. Conversely, if a tight cut
of $G$ is a cut of $G'$ then it is also tight in $G'$.
\end{lem}

%Chen et al.\cite{Lu} show the following result about barriers of a matching covered graph.

\begin{lem}\label{BB}
Let $B$ be a nontrivial barrier of a matching covered graph $G$, and $G[Y]$
be an odd component of $G-B$. Suppose that $G'= G/(\overline{Y}\rightarrow\overline{y})$, and $B'$ is a barrier of $G'$. If
$\overline{y}\in B'$, then $B\cup (B'-\overline{y})$ is a barrier of $G$, and every component of $G'-B'$ is also a component of $G-(B\cup (B'- \overline{y}))$, otherwise $B'$ is also a barrier of $G$.
\end{lem}

\proof Let $H_1,H_2,\ldots,H_{|B|-1},Y$ be all odd  components of $G-B$, and $O_1,O_2,\ldots,O_{|B'|}$ be all odd  components of $G'-B'$. Since $Y$ is an odd  component of $G-B$, every edge that is incident with $\overline{y}$ in $G'$ is incident with one vertex of $B$ in $G$. If $\overline{y}\in B'$, then $B\cup (B'-\overline{y})$ is a barrier of $G$ and $H_1,H_2,\ldots,H_{|B|-1},O_1,O_2,\ldots,O_{|B'|}$ are all odd  components of $G-(B\cup (B'-\overline{y}))$. If $\overline{y}\notin B'$, adjust notation so that $\overline{y}\in O_1$. By Lemmas \ref{2-c} and \ref{C-t}, $G[\overline{Y}]$ is connected. Then $B'$ is also a barrier of $G$ and $G[V(O_1-\overline{y})\cup \overline{Y}],O_2,\ldots,O_{|B'|}$ are all odd  components of $G-B'$ since
\endproof

\begin{lem}\label{B-2}
Let $\{u,v\}$ be a $2$-separation of  matching covered graph $G$ which gives rise to $2$-separation cut $D:=\partial(Y)$. Adjust notation so that $u\in Y$ and $v\in \overline{Y}$. Suppose that $G'= G/(\overline{Y}\rightarrow\overline{y})$, and $B$ is a barrier of $G'$. If
$\overline{y}\in B$, then $(B-\overline{y})+v$ is a barrier of $G$. If
$\overline{y}\notin B$, then $B$ is a barrier of $G$.
\end{lem}

\proof If $\overline{y}\notin B$, then it is obvious that $B$ is a barrier of $G$. If $\overline{y}\in B$,  then $u\notin B$ since $\overline{y}u\in E(G')$ and $G'$ is matching covered graph. Let $L$ be a component of $G'-B$ that contains vertex $u$.
We conclude that  except $\overline{y}u$,  every edge that is incident with $\overline{y}$ in $G'$ is incident with $v$ in $G$. Hence every component of $G'-B$ is also a component of $G-((B-\overline{y})+v)$ except $L$. It is easy to see that $G[V(L)\cup (\overline{Y}-v)]$ is an odd component of $G-((B-\overline{y})+v)$. Hence $(B-\overline{y})+v$ is also a barrier of $G$.\endproof

\section{Proof of the main result}
In this section, the proof of Lemma \ref{O-component} uses arguments similar to that of Theorem 4.1 of \cite{ot}. The statement of Lemma \ref{zy} is the same as Theorem 1.11 of \cite{Lu}.

\begin{lem}\label{O-component}
Let $C:=\partial(X)$ be a nontrivial tight cut of a matching covered graph $G$. If there is  an edge $e:=uv\in C$ such that $u\in X$, $v\in \overline{X}$ and both $G[X-u]$ and $G[\overline{X}-v]$
are connected, then $G$ has a nontrivial barrier that is a proper subset of  $X$ or $\overline{X}$, or has a nontrivial $2$-separation cut that  does not cross $C$.
%Specifically, if $G[\overline{X}+u]$ is 2-connected, then the $2$-separation cut arises from $2$-separation $\{u,z\}$, where \blue{$z\in \overline{X}$.}}
\end{lem}

\proof
%Every nontrivial tight cut in a bipartite graph is a barrier cut. \red{Now suppose that $G$ is nonbipartite.}
By  Lemmas \ref{2-c} and \ref{C-t}, $G/(X\rightarrow x)$ and $G/(\overline{X}\rightarrow \overline{x})$ are 2-connected, whence both $G[X]$ and $G[\overline{X}]$ are connected. The analysis is divided into the following two cases.

%\red{The} analysis may be divided into the following two cases.

%\textbf{
%\red{Case 1:} There is  an edge $e:=uv\in C$ such that $u\in X$, $v\in \overline{X}$ and both $G[X\setminus \{u\}]$ and $G[\overline{X}\setminus \{v\}]$
%are connected.}

\emph{Case 1: $v$ is the only neighbor of $u$ in $\overline{X}$ and $u$ is the only neighbor of $v$
in $X$.}

Let $G':=G-u-v$. Since $G$ is a matching covered graph, $G$ has a
perfect matching $M$ that contains edge $e$. Then $M-e$ is a perfect
matching of $G'$, and we deduce that $G'$ has perfect matchings. Moreover,
$C-e=\partial_{G'}(X-u)=\partial_{G'}(\overline{X}-v)$ is a cut of $G'$ and no edge of $C-e$ is admissible in $G'$.

By Lemma \ref{DM-B}, $G'$ has a DM-barrier
$B'$ such that $B'$, as well as the vertex sets of all the odd components $O_1, O_2, \ldots, O_{|B'|}$
of $G'-B'$, are subsets of one of $X-u$ and $\overline{X}-v$. Adjust notation so that
$B'\subsetneq X-u$ and $O_{i}\subsetneq X-u$ for $i\in \{1,2,\cdots,|B'|\}$.  Since $G[\overline{X}-v]$ is connected, $\overline{X}-v$ is a subset of the vertex set of an even component, say
$L$, of the graph $G'-B'$.

Let $B:= B'\cup\{u\}$. By hypothesis, $u$ is the only vertex of
$X$ adjacent to $v$. Thus, $O_1, O_2, \ldots, O_{|B'|}$ and $G[V(L)+v]$ are also
odd components of $G-B$, implying that $B$ is a nontrivial barrier of $G$ and $B\subsetneq X$.
%Since $C$ is nontrivial and $\overline{X}\subseteq V(L)\cup\{v\}$, $\partial(V(L)\cup\{v\})$ is a nontrivial barrier cut and does not cross $C$.

\emph{Case 2: $u$ has two or more neighbours in $\overline{X}$,
or $v$ has two or more neighbours in $X$.}

Without loss of generality, we assume that $u$ has two or more neighbours in $\overline{X}$. Let $R:=\partial(u)\backslash C$. Now consider the graph
$G'':=G-R$, together with the cut $D:=\partial(X-u)$ in $G''$. Since both $G[X-u]$ and $G[\overline{X}]$ are  connected, the graphs $G''[X-u]$ and $G''[\overline{X}+u]$ are both connected.

Every perfect matching of $G$  containing edge $e$ is also a perfect matching
of $G''$. Thus, $G''$ has perfect matchings. Since $u$ is not adjacent to any vertex of $X$ in $G''$ and $C$ is a tight cut of $G$, no edge of $D$ is admissible in $G''$. By Lemma \ref{DM-B}, $G''$ has a DM-barrier
$B''$ such that $B''$, as well as the vertex sets of all the odd components $H_1, H_2, \ldots, H_{|B''|}$
of $G''-B''$, are subsets of one of $X-u$ and $\overline{X}+u$.

If $V(H_1), V(H_2), \ldots, V(H_{|B''|})$ and $B''$ are proper subsets of $X-u$, then  $\overline{X}+u$  is a subset of the vertex set of an even component, say
$L_1$, of the graph $G''-B''$ since $G[\overline{X}+u]$ is connected. Hence $B:=B''\cup\{u\}$ is a nontrivial barrier of $G$, and $H_1, H_2, \ldots, H_{|B''|}$ and $G[V(L_1)-u]$ are all odd components of $G-B$. As $B''\subsetneq X-u$ and $u\in X$, $B\subsetneq X$.
%Since $C$ is nontrivial and $\overline{X}\subseteq V(L_1)-u$, $\partial(V(L_1)-u)$ is a nontrivial barrier cut that does not cross $C$, and $B\subsetneq X$.

If $V(H_1), V(H_2), \ldots, V(H_{|B''|})$ and $B''$ are proper subsets of  $\overline{X}+u$, then $X-u$ is a subset of the vertex set of an even component, say
$L_2$, of the graph $G''-B''$ since $G[X-u]$ is connected. Then $u$ does not lie in $B''$. Otherwise, $G''-B'' = G-B''$ and $B''$ is a barrier of $G$. Further, $L_2$ is also an even component of $G-B''$, contradicting the condition that $G$ is matching covered. Hence $u$ is a vertex of an even component of $G''-B''$ or $u$ is a vertex of an odd component of $G''-B''$.

\emph{Subcase 2.1: $u$ is a vertex of an even component of $G''-B''$, say $L_3$.}

 Then $L_2=L_3$. Otherwise, $H_1, H_2, \ldots, H_{|B''|}$ and $G[V(L_3)-u]$ are all odd components of $G-(B''\cup\{u\})$, implying that $B''\cup\{u\}$ is a nontrivial barrier of $G$. $L_2$ is an even component of $G''-B''$, implying that $L_2$ is also an even component of $G-(B''\cup\{u\})$, contradicting the condition that $G$ is  matching covered. Hence $L_2=L_3$.  Since $u$  is not adjacent to any vertex of $X-u$ in $G''$ and $|V(L_2)|$ is even, $V(L_2)\cap \overline{X}\neq \emptyset$. For any $w\in V(L_2)\cap \overline{X}$, $B:=B''\cup\{w\}$ is a nontrivial barrier of $G$, and $H_1, H_2, \ldots, H_{|B''|}$ and $G[ V(L_2)-w]$ are all odd components of $G-B$. As $B''\subsetneq \overline{X}$ and $w\in \overline{X}$, $B\subsetneq \overline{X}$.
 %Since $C$ is nontrivial and $X \subseteq V(L_2)-w$, $\partial(V(L_2)-w)$ is a nontrivial barrier cut and does not cross $C$ and $B\subsetneq \overline{X}$.

\emph{Subcase 2.2: $u$ is a vertex of an odd component of $G''-B''$.}

Adjust notation so that $u\in V(H_1)$. Note that $u$ is the only vertex in an odd component of $G''-B''$ that is
adjacent to vertices of $X$ in $G$. Then the barrier $B''$ of $G''$ is also a barrier of $G$, and  $H_2, \ldots, H_{|B''|}$ and $G[V(H_1)\cup V(L_2)]$ are all odd components of $G-B''$.  If $B''$ is nontrivial, we are done with this proof.

If $B''$ is trivial, $u$ lies in $V(H_1)$ and has at least two neighbours in $G''$. Thus there is  also at least
one neighbour of $u$ lying in $V(H_1)$, whence $|V(H_1)|\geq2$. As $H_1$ is an odd component of $G''-B''$,  $|V(H_1)-u|$ is even. Let $z$ denote
the only vertex of $B''$.  We have $\{u,z\}$ is a $2$-separation of $G$, and  $\partial((V(H_1)-u)+z)$ is a $2$-separation cut of $G$. As $(V(H_1)-u)+z\subseteq \overline{X}$, $\partial((V(H_1)-u)+z)$  does not cross $C$. \endproof

%Specifically, if $G[\overline{X}+u]$ is 2-connected, then
%$G[\overline{X}+u-z]$ is connected. As $G''[\overline{X}+u-z]=G[\overline{X}+u-z]$,  $G''[\overline{X}+u-z]$ is connected and %$H_1=G''[\overline{X}+u-z]$, implying that $\{u,z\}$ is a $2$-separation of $G$, and  $\partial(\overline{X})$ is a $2$-separation cut of $G$, i.e. %$C$ is a $2$-separation cut of $G$. \endproof
\vspace{0.3cm}

\noindent {\bf Remark~1}: In the proof of Subcase 2.2 of Lemma \ref{O-component}, if $B''=\{z\}$ and $G[\overline{X}+u]$ is 2-connected, then
$G[\overline{X}+u-z]$ is connected. As $G''[\overline{X}+u-z]=G[\overline{X}+u-z]$,  $G''[\overline{X}+u-z]$ is connected and $H_1=G''[\overline{X}+u-z]$, implying that $\{u,z\}$ is a $2$-separation of $G$, and  $\partial(\overline{X})$ is a $2$-separation cut of $G$, i.e. $C$ is a $2$-separation cut of $G$. Hence under the condition of Lemma \ref{O-component}, if $G[\overline{X}+u]$ is 2-connected, then $G$  has a nontrivial barrier that is a proper subset of  $X$ or $\overline{X}$, or $C$ is $2$-separation cut of $G$ that arises from $2$-separation $\{u,z\}$, where $z\in \overline{X}$.

%By Lemma \ref{O-component} and Remark \ref{rem}, it is not difficult to obtain the following  corollary.
%
%\begin{cor}\label{cor}
%Let $C:=\partial(X)$ be a nontrivial tight cut of matching covered graph $G$. If there is an edge $e:=uv\in C$ such that $u\in X$, $v\in \overline{X}$ and both $G[X-u]$ and $G[\overline{X}-v]$
%are connected, and $G[\overline{X}+u]$ is 2-connected, then $G$  has a nontrivial barrier $B$ of $G$  such that $B\subsetneq X$ or $B\subsetneq \overline{X}$, or $C$ is $2$-separation cut of $G$ that arise in $2$-separation $\{u,x\}$, where $x\in \overline{X}$.
%\end{cor}

%
%\begin{lem}\label{O-component-1}
%Let $C:=\partial(X)$ be a nontrivial tight cut of matching covered graph $G$. If there is no edge $e:=uv\in C$, with $u\in X$ and $v\in \overline{X}$, such that both  $G[X-u]$ and $G[\overline{X}-v]$ are connected, then $G$ has a nontrivial $2$-separation cut that  does not cross $C$, or has a nontrivial barrier $B$ of $G$  such that $B\subsetneq X$ or $B\subsetneq \overline{X}$.
%\end{lem}

\begin{lem}\label{zy}
Let $C:=\partial(X)$ be a nontrivial tight cut of a matching covered graph $G$. Then $G$ has a nontrivial barrier $B$ of $G$  that is a proper subset of $X$ or $\overline{X}$, or has a nontrivial $2$-separation cut that  does not cross $C$.
\end{lem}

\proof If there is an edge $e:=uv\in C$ such that $u\in X$, $v\in \overline{X}$ and both $G[X-u]$ and $G[\overline{X}-v]$
are connected, then by Lemma \ref{O-component},  the result holds.

Now suppose that there is no edge $e:=uv\in C$, with $u\in X$ and $v\in \overline{X}$, such that both  $G[X-u]$ and $G[\overline{X}-v]$ are connected.  Then for any edge $e:=uv\in C$, with $u\in X$ and $v\in \overline{X}$,  $u$ is a cut vertex of $G[X]$ or $v$ is a cut vertex of $G[\overline{X}]$. If $G[X]$ has no cut vertex, then  every vertex of $G[\overline{X}]$ that is incident with an edge of $C$ is a cut vertex of $G[\overline{X}]$ by hypothesis, implying that some cut vertex of $G[\overline{X}]$ is also a cut vertex of $G$. This contradicts that $G$ is matching covered graph. Hence $G[X]$  has two or more blocks. Likewise, $G[\overline{X}]$ has two or more blocks.

\begin{figure}[h]
\centering
\input{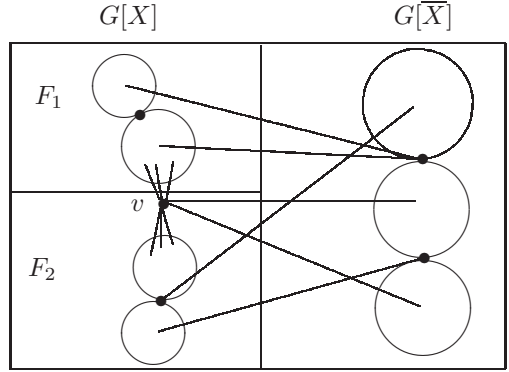}
\caption{The component $F_1$ of $G[X]-v$}
\label{f1}
\end{figure}

Since $G$ is  finite and 2-connected, there is a cut vertex $v$ of $G[X]$ which is incident with an edge of $C$,  such that one component $F_1$ of $G[X]-v$ has no cut vertex of $G[X]$ that is incident with an edge of $C$. Let $F_2:=G[X\setminus V(F_1)]$, as shown in Figure \ref{f1}.
 %Note that for any $e:=v_1v_2\in C$, if  $v_1\in V(F_1)$, then $v_2$ is a cut vertex of $G[\overline{X}]$.

As $v$ is a cut vertex of $G[X]$, $G/(\overline{X}\rightarrow \overline{x})-\{\overline{x},v\}$ is not connected. If $o(G/(\overline{X}\rightarrow \overline{x})-\{\overline{x},v\})>0$, then $o(G/(\overline{X}\rightarrow \overline{x})-\{\overline{x},v\})=|\{\overline{x},v\}|=2$ since $G/(\overline{X}\rightarrow \overline{x})$ is matching covered graph, i.e., $\{\overline{x},v\}$ is a barrier of $G/(\overline{X}\rightarrow \overline{x})$ and $\overline{x}v\notin E(G/(\overline{X}\rightarrow \overline{x}))$. This contradicts that $v$ is incident with an edge of $C$. Hence $o(G/(\overline{X}\rightarrow \overline{x})-\{\overline{x},v\})=0$, implying that $\{\overline{x},v\}$ is a $2$-separation
of $G/(\overline{X}\rightarrow \overline{x})$.

\begin{figure}[htbp]
\centering
\input{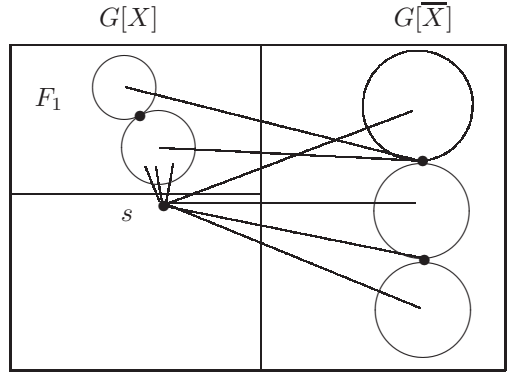}
\caption{Graph $G_2:=G/(V(F_2)\rightarrow s)$}
\label{f2}
\end{figure}

Hence $\partial(V(F_2))$ is a $2$-separation cut of $G/(\overline{X}\rightarrow \overline{x})$.
By Lemma \ref{C-t}, $\partial(V(F_2))$ is also  a tight cut of $G$. Let $G_2:=G/(V(F_2)\rightarrow s)$, as shown in Figure \ref{f2}. By Lemma \ref{C-t}, $C$ is a nontrivial tight cut of $G_2$ and $V(F_1)+s$ and $\overline{X}$ are two shores of $C$.

If $G_2[\overline{X}+s]$ has two or more blocks, there is a block $K$ that satisfies $s\notin V(K)$ and contains precisely one cut vertex, say $k$, of $G_2[\overline{X}+s]$. Note that $s$ is not cut vertex of $G_2[\overline{X}+s]$ since $G_2[\overline{X}]$ is connected.
By the definition of $F_1$, every edge of $C$ that is incident with a vertex of $F_1$ is incident with a cut vertex of $G_2[\overline{X}]$, implying that $k$ is also a cut vertex of $G_2$. This contradicts that $G_2$ is matching covered graph. Hence $G_2[\overline{X}+s]$ is 2-connected.

Since $G_2[\overline{X}+s]$ is 2-connected, there is a vertex $w$ of $\overline{X}$ which is not a cut vertex of $G[\overline{X}]$ and adjoins $s$. It is easy to see that $sw\in C$ and both $G_2[\overline{X}-w]$ and $F_1$ are  connected. By Remark~1,  $G_2$  has a nontrivial barrier $B$ satisfying $B\subsetneq \overline{X}$ or $B\subsetneq V(F_1)+s$, or $C$ is $2$-separation cut of $G_2$ that arises from $2$-separation $\{s,z\}$, where $z\in \overline{X}$. Next, two cases should be considered.

\emph{Case 1: There is a nontrivial barrier $B$ in $G_2$, where $B\subsetneq \overline{X}$ or $B\subsetneq V(F_1)+s$.}

If $s\notin B$,  then $B$  is also a barrier of $G$. If $s\in B$, then $B\subsetneq V(F_1)+s$,  implying that $\overline{X}$ is a
subset of the vertex set of an odd component, say $L$, of $G_2-B$. Since $v$ is a cut vertex of $G[X]$, every edge of  $G_2[V(F_1)+s]$ that is incident with $s$ is  incident with $v$ in $G$. Hence $B':=(B-s)+v$ is  a barrier of $G$, $G[V(L)\cup V(F_2-v)]$ is an odd component of $G-B'$ and all odd components of $G_2-B$ that are subgraphs of $F_1$ are also odd components of $G-B'$. Obviously, $B'\subsetneq X$.

\emph{Case 2: $C$ is $2$-separation cut of $G_2$ that arise from $2$-separation $\{s,z\}$, where $z\in \overline{X}$}

Since  every edge that joins a vertex of $F_1$ to  $s$ is  incident with $v$ in $G$, we have $\{v,z\}$ is a $2$-separation of $G$, and $\partial(V(F_1)+v)$ is a 2-separation cut of $G$. Since $V(F_1)+v \subsetneq X$, $\partial(V(F_1)+v)$  does not cross $C$.\endproof

%The following result follows directly from Lemmas \ref{O-component} and \ref{O-component-1}.

%\begin{cor}\label{zy}
%Let $C:=\partial(X)$ be a nontrivial tight cut of matching covered graph $G$ \red{which is} not an ELP-cut. Then $G$ has a nontrivial $2$-separation cut that  does not cross $C$, or has a nontrivial barrier $B$ of $G$  such that $B\subsetneq X$ or $B\subsetneq \overline{X}$.
%\end{cor}

\vspace{0.3cm}
Now, we are going to prove our main result, Theorem \ref{2-S}.

\textbf{Proof of Theorem \ref{2-S} } Let $C:=\partial(X)$ be a tight cut in a matching covered graph $G$ which is not an ELP-cut. We first prove the following two claims.

\emph{ \textbf{Claim 1}: If $G$ has no barrier that is a proper subset of $X$ or $\overline{X}$, then
there is a sequence $G_1=G, G_2,\ldots,G_r, r\geq2$  of matching covered graphs,
such that for $i=1, 2,\ldots, r-1$, $G_i$ has a $2$-separation cut $C_i$,  $G_{i+1}$ is a
$C_i$-contraction of $G_i$, and  $C$ is a $2$-separation cut of $G_r$.}

 By Lemma \ref{zy}, $G$ has a $2$-separation cut $C_1$ that does not cross $C$.  Let $G_{2}$ be the
$C_1$-contraction of $G_1=G$ that contains $C$. By Lemma \ref{B-2},  $G_{2}$ has no nontrivial barrier that is a proper subset of one of two shores of $C$ in $G_{2}$.
By Lemma \ref{zy}, there is a $2$-separation cut $C_{2}$ that  does not cross $C$ in $G_{2}$. Applying Lemmas \ref{B-2} and \ref{zy} repeatedly, we can recursively get $G_{i+1}$ is a
$C_i$-contraction of $G_i$ for $i=1,\ldots$, where $C_i$ is a $2$-separation cut and  does not cross $C$.

Since $V(G)$  is finite and $C$ is not an ELP-cut in $G$, there is an index $r$ such that $C$ is a $2$-separation cut of $G_r$, where $r\geq 2$.  Claim 1 holds.

\vspace{0.3cm}

\emph{\textbf{Claim 2}: If $G$ has a nontrivial barrier that is a proper subset of $\overline{X}$(resp. $X$), then $G$ has a nontrivial barrier cut $C^*$ such that one of its shores, say $Y$, is a superset of $X$ and $G':=G/\overline{Y}$ has no nontrivial barrier that is a proper subset of $V(G')\setminus X$(resp. $V(G')\setminus \overline{X}$).}

Let $B$ be a nontrivial barrier of $G$. Without loss of generality, we assume that $B\subsetneq \overline{X}$. Since $G[X]$ is connected, there is an odd component $G[Y]$ of $G-B$ such that $X\subseteq Y$. Choose $B$ so that $|Y|$ is as small as possible. Let $G':=G/\overline{Y}$. Then $G'$ has no nontrivial barrier that is a proper subset of $V(G')\setminus X$. Otherwise, we assume that $B'\subsetneq V(G')\setminus X$ is a nontrivial barrier of $G'$ and $G[Y']$ is an odd component of $G'-B'$ with  $X\subseteq Y'$. By Lemma \ref{BB}, there is a nontrivial barrier $B''\subsetneq\overline{X}$ of $G$ such  that $G[Y']$ is an odd component of $G-B''$. Obviously, $Y'\subsetneq Y$, in contradiction to the minimality of $Y$. Hence Claim 2 holds.

\vspace{0.3cm}
If $G$ has no barrier that is a proper subset of $X$ or $\overline{X}$, the assertion holds by Claim~1. If $G$ has a nontrivial barrier $B$ that is a proper subset of $X$ or $\overline{X}$, by symmetry, we may assume $B\subsetneq \overline{X}$. By Claim~2, then $G$ has a nontrivial barrier cut $C^*:=\partial(Y)$ and $X\subsetneq Y$,  and $G':=G/\overline{Y}$ has no nontrivial barrier that is a proper subset of $V(G')\setminus X$. Next if $G'$ also has no nontrivial barrier that is a proper subset of $X$, then the assertion holds by Claim~1. If $G'$ has a nontrivial barrier that is a proper subset of $X$, by Claim~2 and Lemma \ref{BB}, then $G'$ has a nontrivial barrier cut $C':=\partial(Z)$ and $V(G')\setminus X \subsetneq Z$,  and $G'':=G'/\overline{Z}$ has no nontrivial barrier that is a proper subset of one of the two shores of $C$ in $G''$. By Claim~1, the assertion holds.\endproof

\acknowledgements
\label{sec:ack}
We would like to thank the responsible editor and the anonymous referees for their the constructive comments and kind suggestions
on improving the representation of the paper.

\nocite{*}
\bibliographystyle{abbrvnat}
% use the following instead if you encounter problems
%\bibliographystyle{alpha}
\bibliography{A-note-on-tight-cuts-in-matching-covered-graphs}

\begin{thebibliography}{8}
\providecommand{\natexlab}[1]{#1}
\providecommand{\url}[1]{\texttt{#1}}
\expandafter\ifx\csname urlstyle\endcsname\relax
  \providecommand{\doi}[1]{doi: #1}\else
  \providecommand{\doi}{doi: \begingroup \urlstyle{rm}\Url}\fi

\bibitem[Bondy and Murty(2008)]{Bondy}
J.~Bondy and U.~S.~R. Murty.
\newblock \emph{Graph Theory}.
\newblock Springer, New York, 2008.

\bibitem[Carvalho et~al.(2002)Carvalho, Lucchesi, and Murty]{G-e}
M.~H. Carvalho, C.~L. Lucchesi, and U.~S.~R. Murty.
\newblock On a conjecture of lov\'asz concerning bricks. ii. bricks of finite
  characteristic.
\newblock \emph{J. Combin. Theory Ser. B}, 85\penalty0 (1):\penalty0 137--180,
  2002.

\bibitem[Carvalho et~al.(2018)Carvalho, Lucchesi, and Murty]{ot}
M.~H. Carvalho, C.~L. Lucchesi, and U.~S.~R. Murty.
\newblock On tight cuts in matching covered graphs.
\newblock \emph{J. Comb.}, 9\penalty0 (1):\penalty0 163--184, 2018.

\bibitem[Chen et~al.(2020)Chen, Feng, Lu, Lucchesi, and Zhang]{Lu}
G.~T. Chen, X.~Feng, F.~L. Lu, C.~L. Lucchesi, and L.~Z. Zhang.
\newblock Laminar tight cuts in matching covered graphs.
\newblock preprint, 2020.
\newblock available at
  \href{https://arxiv.org/abs/2003.08622}{\texttt{https://arxiv.org/abs/2003.08622}}.

\bibitem[Edmonds et~al.(1982)Edmonds, Lov\'asz, and Pulleyblank]{Lov}
J.~Edmonds, L.~Lov\'asz, and W.~R. Pulleyblank.
\newblock Brick decompositions and the matching rank of graphs.
\newblock \emph{Combinatorica}, 2\penalty0 (2):\penalty0 247--274, 1982.

\bibitem[Lov\'asz(1987)]{ML}
L.~Lov\'asz.
\newblock Matching structure and the matching lattice.
\newblock \emph{J. Combin. Theory Ser. B}, 43:\penalty0 187--222, 1987.

\bibitem[Plummer(1980)]{2}
M.~D. Plummer.
\newblock On n-extendable graphs.
\newblock \emph{Discrete Mathematics}, 31:\penalty0 201--210, 1980.

\bibitem[Szigeti(2002)]{S}
Z.~Szigeti.
\newblock Perfect matchings versus odd cuts.
\newblock \emph{Combinatorica}, 22:\penalty0 575--589, 2002.

\end{thebibliography}
\label{sec:biblio}

\end{document}